\documentclass[a4paper, 11pt]{article}

\usepackage[latin1]{inputenc}
\usepackage{makeidx }
\usepackage[french]{babel}
\usepackage{amsmath,amsthm,amsfonts,amssymb}
\usepackage{newlfont}
\usepackage{showidx }
\usepackage[dvips]{graphicx}
\newtheorem{theo}{Théorème}

\def\g{\mathfrak{g}}

\def\tr{\mathrm{tr}}
\def\d{\mathrm{d}}
\def\ad{\mathrm{ad}}

\begin{document}
\bibliographystyle{alpha}
\author{Luc Albert
\footnote{Lycée Massena,  2 avenue Félix Faure,
06050 Nice Cedex 1, luc@albert1.net }
\and Pascale Harinck\footnote{ CM Laurent Schwartz,
Ecole polytechnique, CNRS,
91128 Palaiseau Cédex
FRANCE, harinck@math.polytechnique.fr}
\and Charles Torossian\footnote{Institut Mathématiques
de Jussieu, CNRS , Université Paris 7, Site Chevaleret, Case
7012, 2 place Jussieu,  75205 Paris Cedex 13, FRANCE, Email: torossian@math.jussieu.fr}}
\title{Solution non universelle pour le problème KV-78}

\maketitle

{\noindent \bf{R\'esum\'e:}} En 1978, M. Kashiwara et M. Vergne \cite{KV} ont conjecturé qu'il était possible d'écrire d'une certaine façon la formule de Campbell-Hausdorff (conjecture KV-78) pour qu'une formule de trace soit vérifiée. Ils proposaient dans le même article une solution à ce problème dans le cas résoluble. La conjecture KV-78 a été démontrée par Alekseev-Meinrenken \cite{AM} en utilisant les travaux du troisième  auteur. On renvoie à \cite{To2} pour un panorama sur ces questions. On montre dans cette note que la solution du cas résoluble proposée dans \cite{KV} ne convient pas en degré 8 pour les algèbres de Lie générales. Nos calculs ont été fait par ordinateur. Les programmes du premier (L.A.) et troisième auteur (C.T.) ont permis par ailleurs de vérifier jusqu'à l'ordre 16 la conjecture de \cite{AT} concernant l'égalité entre les algèbres de Lie de Drinfeld $\mathfrak{grt}_1$ et $\widehat{kv_2}$ définie dans \cite{AT}.\\

\noindent \textbf{Abstract :} In 78' M. Kashiwara and Vergne conjectured some property on the Campbell-Hausdorff series in such way a trace formula is satisfied. They proposed an explicit solution in the case of solvable Lie algebras. In this note we prove that this \textit{solvable solution} is not universal. Our method is based on computer calculation. Furthermore our programs prove up to degree 16, Drinfeld's Lie algebra  $\mathfrak{grt}_1$ coincides with the Lie algebra $\widehat{kv_2}$ defined in \cite{AT}.\\

 \section{Notations-Rappels}

 Soit $\g$ une alg\`ebre de Lie de dimension finie sur  $\mathbf{R}$.
D'apr\`es le  th\'eor\`eme de Lie, il existe un groupe de Lie r\'eel  $G$,
connexe et simplement connexe d'alg\`ebre de Lie $\g$ et  une application exponentielle
not\'ee $\exp_{\g}$ qui d\'efinit un diff\'eomorphisme local en $0\in \g$
sur $G$.

Pour $X, Y$
proches de $0$ dans $\g$, il existe une s\'erie en des
polyn\^omes de Lie,  convergente et \`a valeurs dans~$\g$, not\'ee $Z(X, Y)$ telle que
l'on ait
$$\exp_{\g}(X)\cdot_G\exp_{\g}(Y)=\exp_{\g}\big( Z(X,Y)\big).$$
C'est la série de Campbell-Hausdorff.
Les premiers termes  sont bien connus et s'\'ecrivent
\begin{multline}
Z(X,Y)=X+Y +\frac{1}{2} [X,Y] + \frac{1}{12}[X,[X,Y]]
+\frac{1}{12}[Y,[Y,X]]+\\
\frac{1}{48}[Y,[X,[Y,X]]]-\frac{1}{48}[X,[Y,[X,Y]]]+ \cdots.
\end{multline}

Cette formule est universelle, c'est donc aussi une formule dans l'algèbre de Lie libre engendrée par $X, Y$ que l'on complète par le degré des crochets.

La conjecture combinatoire de Kashiwara--Vergne \cite{KV}  a été démontrée par Alekseev-Meinrenken  \cite{AM} en utilisant les résultats du troisième auteur \cite{To} et s'\'enonce de
la mani\`ere suivante.

\begin{theo}[Conjecture KV-78]\label{conjKV}Notons $Z(X, Y)$ la s\'erie de Campbell-Hausdorff. Il existe des s\'eries $F(X,Y)$ et $G(X,Y)$ sur $\g\oplus\g$ sans
terme constant et \`a valeurs dans $\g$ telles que l'on ait
\begin{equation}\label{KV1}
X+Y-\log(\exp_\g(Y)\cdot \exp_\g (X)) =\big(1-e^{-\ad X}\big)F(X,Y)+ \big(e^{\ad Y}
-1\big)G(X,Y)
\end{equation}
et telle que l'identit\'e de trace suivante soit v\'erifi\'ee
\begin{equation}\label{KV2}\tr_{\g}\left(\ad X\circ \partial_X F+\ad Y\circ \partial_Y G\right)=
\frac{1}{2}\tr_{\g}\left( \frac{\ad X}{e^{\ad
X}-1}+\frac{\ad Y}{e^{\ad Y}-1} -\frac{\ad Z(X, Y)}{e^{\ad
Z(X, Y)}-1}-1\right).
\end{equation}
\end{theo}

Si $(F, G)$ est un couple solution des équations de KV-78 alors  $$\left(G(-Y, -X), F(-Y, -X)\right)$$
est aussi une solution. On peut donc rechercher des solutions sym\'etriques,
c'est-\`a-dire v\'erifiant $G(X, Y)=F(-Y, -X)$.
\section{La solution \textit{résoluble} de Kashiwara--Vergne}

Dans leur article, Kashiwara et Vergne proposent, dans le cas \textit{résoluble},  une solution sym\'etrique  $(F^0, G^0)$ de s\'eries de Lie universelles. Il existe dans le cas quadratique,  deux solutions \cite{Ve, AM2} qui sont toutefois différentes de la solution résoluble. Dans \cite{bur}, est proposée une solution pour la première équation, qu'il serait intéressant de comparer avec la solution résoluble.\\

 Suivant  \cite{Rou86} le couple de solutions proposé dans  \cite{KV} se décrit de la manière suivante: notons  $\psi$ la  fonction analytique au voisinage de $0$  d\'efinie par  $$\psi(z)= \frac{e^z-1-z}{(e^z-1)(1-e^{-z})}.$$
Soit $Z(t)=Z(tX, tY)$~; posons  $$F^1(X, Y)= \left(\int_0^1 \frac{1-e^{-t \,\ad X}}{1-e^{-\ad
X}} \circ \psi(\ad Z(t)) \d t\right) (X+Y)$$ et $G^1(X, Y)=F^1(-Y, -X)$. Posons alors
\begin{equation}F^0(X,
Y)=\frac12\left( F^1(X, Y) +e^{\ad X} F^1(-X, -Y)\right) +\frac14\left( Z(X, Y) -X\right)\end{equation}\label{Fzero} et
$G^0(X, Y)=F^0(-Y, -X)$.\\

Par construction (\cite{KV}, \cite{Rou86}) ce couple $(F^0, G^0)$ v\'erifie l'\'equation~(\ref{KV1}) pour toute algèbre de Lie (pas uniquement les algèbres résolubles). Le couple $(F^0, G^0)$ v\'erifie aussi l'\'equation de trace~(\ref{KV2}) dans le cas résoluble. On appellera ce couple:  \textit{solution résoluble}.\\

On va montrer que cette solution ne convient pas pour les algèbres de Lie générales à partir du degré 8.\\

Les premiers termes de la série ont été calculés, à la main, jusqu'à l'ordre~4  par Rouvière dans \cite{Rou86}, mais Rouvière a mené les calculs jusqu'à l'ordre 6. On trouve
\begin{multline}
F^0(X, Y)=
\frac14 Y + \frac1{24}[X,Y]
-\frac1{48}[X, [X,Y]]-\frac1{48}[Y,[X,Y]]\\-\frac1{180}[X, [X,[X,Y]]]-\frac1{480}[Y,[X,[X,Y]]]+\frac1{360}[Y,[Y, [X,Y]]] + \ldots
\end{multline}

Comme notre couple est symétrique on a $G^0(X,Y)=F(-Y, -X)$.\\
\section{Présentation des programmes}

\subsection{Procédures Maple et calcul en Matlab}

Notre programme en Maple, écrit par les deux derniers auteurs (P. Harinck et C.~Torossian) est disponible sur : \\  \verb"http://www.institut.math.jussieu.fr/~torossian/publication/pubindex.html"

Notre programme calcule pour $N$ donné les termes d'ordre $N$ de $F^0(X,Y)$  grâce à la procédure  \verb"procedure KV(N)".\\

Ces termes sont obtenus sous formes de crochets itérés. Afin d'éviter les redondances et pour économiser un peu de place, nous avons présentés, dans cette note, les résultats dans la base de Hall. Les calculs ont été fait sous Matlab, grâce au package \verb"Lafree" dans \verb"Diffman" et disponible sur Internet.  Ce package est très facile d'utilisation, il permet de générer une base de Hall et de calculer les coordonnées de tout mot de Lie écrit dans cette base grâce à la fonction \verb"getvector".\\

Notre programme Maple calcule par la procédure \verb"procedure SM(N)", dans $Ass_2$ l'algèbre  associative libre engendrée par deux variables $1=\ad X, 2=\ad Y$,  les termes d'ordre $N$ dans  $$\frac{1}{2}\left( \frac{\ad X}{e^{\ad
X}-1}+\frac{\ad Y}{e^{\ad Y}-1} -\frac{\ad Z(X, Y)}{e^{\ad
Z(X, Y)}-1}-1\right). $$ On détermine alors les coordonnées dans une base des mots cycliques. La base choisie est indiquée dans la note.

Nous avons remplacé comme dans \cite{AT} la condition de trace, par une égalité dans l'espace vectoriel gradué  $cy_2 : = Ass_2/< ab-ba, a, b \in Ass_2>$, c'est à dire que le mot $x_1 x_2 \ldots  x_n$ est identifié au mot $x_2 x_3 \ldots x_n x_1$.  C'est l'espace vectoriel engendré par les mots cycliques\footnote{Dans \cite{bur} cet espace est aussi considéré.}. C'est en effet la seule relation universelle vérifiée par la trace. On n'a pas cherché à optimiser informatiquement cette partie, car notre objectif était l'ordre d'atteindre l'ordre 8 ou 10.\\

\subsection{Programme Caml et résultats complémentaires dus à L.~Albert et C.~Torossian}

Le premier auteur (L.A.) a écrit un programme \verb"CharlesKV2.ML" en Caml, indépendant du programme Maple, disponible sur la page web du troisième auteur (C.T.)\\  \verb"http://www.institut.math.jussieu.fr/~torossian/publication/pubindex.html".
Le programme  calcule $$\ad X\circ \partial_X F+\ad Y\circ \partial_Y G$$ dans la base cyclique générée par le programme.\\

Ce programme a permis de vérifier jusqu'à l'ordre 16 la conjecture de \cite{AT} concernant l'égalité des dimensions entre l'algèbre de Lie $\widehat{kv_2}$ définie dans \cite{AT} et l'algèbre $\mathfrak{grt}_1$ définie par Drinfeld \cite{Dr}. Plus précisément on calcule pour un degré $N$ donné  la dimension du noyau du système qui définit l'algèbre de Lie $kv_2$ : $(A, B)$ couple de mots de Lie en deux variables $X, Y$ de degré $N$  vérifiant les deux équations suivantes

 \begin{eqnarray}
[X, A(X, Y)] + [Y, B(X, Y)]=0\\ \tr_\g \Big( \partial_X A(X, Y) \circ \ad X  +  \partial_Y B(X, Y) \circ \ad Y \Big)=0.
\end{eqnarray}

Grâce au package \verb"Lafree" on calcule la matrice correspondant à la première équation. Grâce au programme Caml du premier auteur, on peut déterminer la matrice correspondant à la deuxième équation. La complexité croit de manière exponentielle et les programmes sont fortement récursifs. Il en résulte des difficultés de gestion de mémoire et de temps de calculs.\footnote{Pour N=16, il a fallu 18h sur un ordinateur puissant pour obtenir la matrice de la seconde équation. La matrice occupe 100MB de fichier texte.} Sur la page web du troisième auteur (C. T.), on trouvera les matrices du système $kv_2$ pour $N=13,14,15,16$ sous forme \textit{sparse} zippée.\footnote{Sous les conseils de Roman Pearce, le calcul du rang se fait par réduction modulo $p$ à partir de  $N=15$ car les matrices sont de taille trop importante. Cette réduction modulo $p$ ne fait que majorer la dimension du noyau. Par ailleurs les résultats de \cite{AT} et \cite{Ra} (cf. plus loin) minorent la dimension.  Cela permet de conclure sur la dimension exacte.} \\

\noindent En degré $8$, on obtient une unique solution.

\begin{multline}
A(x,y)=   4*[[x,y],[y,[x,[x,[x,[x,y]]]]]] + 6*[[x,y],[y,[y,[x,[x,[x,y]]]]]]\\
      - 6*[[x,y],[y,[y,[y,[x,[x,y]]]]]] - 4*[[x,y],[y,[y,[y,[y,[x,y]]]]]]\\
      - 4*[[x,y],[[x,y],[x,[x,[x,y]]]]] + 6*[[x,[x,y]],[y,[x,[x,[x,y]]]]]\\
      + 6*[[x,[x,y]],[y,[y,[x,[x,y]]]]] - 3*[[x,[x,y]],[y,[y,[y,[x,y]]]]]\\
      + 2*[[x,[x,y]],[[x,y],[x,[x,y]]]] - 3*[[x,[x,y]],[[x,y],[y,[x,y]]]]\\
      + 4*[[y,[x,y]],[x,[x,[x,[x,y]]]]] + 9*[[y,[x,y]],[y,[x,[x,[x,y]]]]]\\
      - 12*[[y,[x,y]],[y,[y,[x,[x,y]]]]] - 10*[[y,[x,y]],[y,[y,[y,[x,y]]]]]\\
      + 9*[[y,[x,y]],[[x,y],[x,[x,y]]]] - 3*[[y,[x,y]],[[x,y],[y,[x,y]]]]\\
      - 6*[[x,[x,[x,y]]],[y,[x,[x,y]]]] - 3*[[x,[x,[x,y]]],[y,[y,[x,y]]]]\\
      + 3*[[y,[x,[x,y]]],[y,[y,[x,y]]]]
      \end{multline}
\begin{multline}
B(x,y)=     -4*[[x,y],[x,[x,[x,[x,[x,y]]]]]] - 6*[[x,y],[y,[x,[x,[x,[x,y]]]]]]\\
      + 6*[[x,y],[y,[y,[x,[x,[x,y]]]]]] + 4*[[x,y],[y,[y,[y,[x,[x,y]]]]]]\\
      - 12*[[x,y],[[x,y],[x,[x,[x,y]]]]] + 18*[[x,y],[[x,y],[y,[x,[x,y]]]]]\\
      + 12*[[x,y],[[x,y],[y,[y,[x,y]]]]] - 10*[[x,[x,y]],[x,[x,[x,[x,y]]]]]\\
      - 12*[[x,[x,y]],[y,[x,[x,[x,y]]]]] + 9*[[x,[x,y]],[y,[y,[x,[x,y]]]]]\\
      + 4*[[x,[x,y]],[y,[y,[y,[x,y]]]]] - 9*[[x,[x,y]],[[x,y],[x,[x,y]]]]\\
      - 3*[[y,[x,y]],[x,[x,[x,[x,y]]]]] + 6*[[y,[x,y]],[y,[x,[x,[x,y]]]]]\\
      + 6*[[y,[x,y]],[y,[y,[x,[x,y]]]]] + 9*[[y,[x,y]],[[x,y],[x,[x,y]]]]\\
      + 4*[[y,[x,y]],[[x,y],[y,[x,y]]]] - 3*[[x,[x,[x,y]]],[y,[x,[x,y]]]]\\
      + 3*[[x,[x,[x,y]]],[y,[y,[x,y]]]] + 6*[[y,[x,[x,y]]],[y,[y,[x,y]]]]
      \end{multline}

On présente ci-dessous notre tableau de résultats (dus à L. Albert et C. Torossian), où on a fait figurer la dimension $d_N$ en degré $N$ de l'algèbre de Lie libre en deux générateurs $Lie_2$, la dimension $c_N$ en degré $N$ de l'espace vectoriel des mots cycliques en deux générateurs $cy_2$. Notre système $kv_2(N)$ en degré $N$ est de taille $((d_{N+1}+ c_N) \times 2d_N)$. Pour $N=16$ la matrice est de taille $(11826, 8160)$.\\

D'après \cite{AT} $E:=\widehat{kv_2}/kv_2$ est concentré en degré impair et on a $\dim E_{2n+1}=1$ pour $n\geq 1$. Les dimensions de $\widehat{kv_2}(N)$ s'en déduisent alors. En degré $N=1$, on trouve un élément  central dans $\widehat{kv_2}$, noté $t=(Y,X)$ dans \cite{AT}.

$$\begin{array}{|c|c|c|c|c|c|c|c|c|c|c|c|c|c|c|c|c|}\hline
  N & 1 &2 & 3 & 4 & 5 & 6 & 7 & 8 & 9&10& 11&12&13&14&15&16\\\hline
  \mathrm{Lie_2(N)} &  2 & 1 & 2 & 3 & 6 & 9 & 18 & 30&56&99&186&335&630&1161&2182&4080\\\hline
   \mathrm{cy_2(N)} &  2 & 3 & 4 & 6 & 8 & 14 & 20 & 36&60&108&188&352&632&1182&2192&4116\\\hline
\mathrm{kv_2(N)}& 1& 0& 0&0&0&0&0&1&0&1&1&2&2&3&3&5\\\hline
 \mathrm{\widehat{kv_2}(N)}& 1& 0& 1&0&1&0&1&1&1&1&2&2&3&3&4&5\\\hline\end{array}$$

Dans \cite{Ra} Racinet  introduit une algèbre graduée $\mathfrak{drm}_0$ dont les dimensions jusqu'en degré $N=19$ sont calculées dans \cite{ENR}. Ces trois auteurs montrent, jusqu'en degré $N=19$, l'inclusion  $\mathfrak{drm}_0 \subset \mathfrak{grt}_1$, où $\mathfrak{grt}_1$ est l'algèbre de Lie introduite par Drinfeld \cite{Dr} . Par ailleurs d'après \cite{AT} on a l'inclusion en tout degré $\mathfrak{grt}_1 \subset \widehat{kv_2}$.

Les dimensions de $\mathfrak{dmr}_0$ sont identiques \cite{ENR}   à celles de $\widehat{kv_2}$ (pour $N\geq 2$), on en déduit l'égalité $\mathfrak{dmr}_0= \mathfrak{grt}_1=\widehat{kv_2}/<t>$ jusqu'en degré $N=16$.

L'algèbre de Drinfeld est conjecturalement égale à une algèbre de Lie libre engendrée par des générateurs en degré impair $\sigma_3, \sigma_5, ....$. L'élément de degré $8$ ci-dessus correspond au terme $[\sigma_3, \sigma_5]$.
\section{Equations de Trace pour le couple $(F^0, G^0)$}
On va comparer,  pour le couple symétrique $(F^0, G^0)$  de l'équation (\ref{Fzero}), les deux membres de  l'équation de trace~(\ref{KV2}) grâce aux deux programmes indépendants.
\subsection{Equation de trace à l'ordre 5 et 7}
Voici les termes d'ordre $5, 7$ pour $F^0(X,Y)$ calculés par le programme Maple.\\

\noindent \textbf{Ordre 5:} Au terme $ \frac1{2880} $ près il vaut
\begin{multline}\nonumber[X,[X,[X,[X,Y]]]] + 8*[Y,[X,[X,[X,Y]]]] + 8*[Y,[Y,[X,[X,Y]]]]\\
      + [Y,[Y,[Y,[X,Y]]]] + 6*[[X,Y],[X,[X,Y]]] + 2*[[X,Y],[Y,[X,Y]]]
\end{multline}

\noindent \textbf{Ordre 7: } Au terme $ \frac1{1209600}$ près il vaut
\nonumber\begin{multline}
-10*[X,[X,[X,[X,[X,[X,Y]]]]]] - 120*[Y,[X,[X,[X,[X,[X,Y]]]]]]\\
      - 380*[Y,[Y,[X,[X,[X,[X,Y]]]]]] - 380*[Y,[Y,[Y,[X,[X,[X,Y]]]]]]\\
      - 120*[Y,[Y,[Y,[Y,[X,[X,Y]]]]]] - 10*[Y,[Y,[Y,[Y,[Y,[X,Y]]]]]]\\      - 150*[[X,Y],[X,[X,[X,[X,Y]]]]] - 940*[[X,Y],[Y,[X,[X,[X,Y]]]]]\\
      - 960*[[X,Y],[Y,[Y,[X,[X,Y]]]]] - 210*[[X,Y],[Y,[Y,[Y,[X,Y]]]]]\\
      - 240*[[X,Y],[[X,Y],[X,[X,Y]]]] - 60*[[X,Y],[[X,Y],[Y,[X,Y]]]]\\
      - 60*[[X,[X,Y]],[X,[X,[X,Y]]]] + 30*[[X,[X,Y]],[Y,[X,[X,Y]]]]\\
      + 360*[[X,[X,Y]],[Y,[Y,[X,Y]]]] - 740*[[Y,[X,Y]],[X,[X,[X,Y]]]]\\
      - 1170*[[Y,[X,Y]],[Y,[X,[X,Y]]]] - 180*[[Y,[X,Y]],[Y,[Y,[X,Y]]]]
\end{multline}

Pour les deux couples d'ordre 5 et 7, le membre de gauche de~(\ref{KV2}) calculé par le programme Caml, donne $0$ dans les bases cycliques. Il est donc égal au  le membre de droite qui est nul aussi pour des raisons de parité. En effet la série de Bernoulli est paire à partir de l'ordre 2, de plus   $Z(-X, -Y)=-Z(Y,X)$ et $Z(Y,X)$ est  conjugué à $Z(X, Y)$. Donc le membre de droite de~(\ref{KV2}) n'a que des termes en degré total pair.

\subsection{Equation de trace à l'ordre 6}

Voici les termes d'ordre $6$ pour $F^0(X,Y)$ calculés par le programme Maple.\\

\noindent \textbf{Ordre 6 :} Au terme $\frac1{120960}$ près il vaut
\small\begin{multline}\nonumber24*[X,[X,[X,[X,[X,Y]]]]] + 69*[Y,[X,[X,[X,[X,Y]]]]] + 20*[Y,[Y,[X,[X,[X,Y]]]]]\\
      - 39*[Y,[Y,[Y,[X,[X,Y]]]]] - 12*[Y,[Y,[Y,[Y,[X,Y]]]]] + 144*[[X,Y],[X,[X,[X,Y]]]]\\
    + 114*[[X,Y],[Y,[X,[X,Y]]]] + 18*[[X,Y],[Y,[Y,[X,Y]]]] - 18*[[X,[X,Y]],[Y,[X,Y]]]
\end{multline}
\normalsize

Considérons la base des mots cyclique de longueur 6 suivante :\\

\verb"C6[1] := [1, 1, 1, 1, 1, 1]; C6[2] := [1, 1, 1, 1, 1, 2];"

 \verb"C6[3] := [1, 1, 2, 1, 1, 2]; C6[4] := [1, 1, 1, 2, 1, 2];"

 \verb"C6[5] := [1, 2, 1, 2, 1, 2]; C6[6] := [1, 1, 1, 1, 2, 2];"

 \verb"C6[7] := [1, 1, 2, 2, 1, 2]; C6[8] := [1, 1, 2, 1, 2, 2];"

  \verb"C6[9] := [1, 2, 2, 1, 2, 2]; C6[10] := [1, 1, 1, 2, 2, 2];"

  \verb"C6[11] := [1, 2, 1, 2, 2, 2]; C6[12] := [1, 1, 2, 2, 2, 2];"

   \verb"C6[13] := [1, 2, 2, 2, 2, 2]; C6[14] := [2, 2, 2, 2, 2, 2]"

\vspace{0,3cm}

Notre programme en Maple  a calculé, au facteur $\frac1{120960}$ près, les coordonnées dans  la base cyclique du terme de droite de~(\ref{KV2}). On trouve :

$$[0,-12,36,-96,-144,30,72,72,36,-40,-96,30,-12,0]$$
Le programme en Caml  calcule les coordonnées dans la base cyclique ci-dessus, le  membre de gauche de~(\ref{KV2}). On trouve au facteur  $\frac1{120960}$ près
$$[0, -12, 36, -96, -144, 30, 72, 72, 36, -40, -96, 30, -12, 0]$$
Il y a donc bien égalité à l'ordre $6$ dans l'équation de trace~(\ref{KV2}).

\subsection{Défaut de l'équation de trace à l'ordre 8}

\textbf{Ordre 8: } Au terme $\frac1{ 21772800}$ près il vaut

\small \begin{multline} \nonumber-144*[X,[X,[X,[X,[X,[X,[X,Y]]]]]]] - 666*[Y,[X,[X,[X,[X,[X,[X,Y]]]]]]]\\
      - 1116*[Y,[Y,[X,[X,[X,[X,[X,Y]]]]]]] - 315*[Y,[Y,[Y,[X,[X,[X,[X,Y]]]]]]]\\
      + 612*[Y,[Y,[Y,[Y,[X,[X,[X,Y]]]]]]] + 414*[Y,[Y,[Y,[Y,[Y,[X,[X,Y]]]]]]]\\
      + 72*[Y,[Y,[Y,[Y,[Y,[Y,[X,Y]]]]]]] - 2160*[[X,Y],[X,[X,[X,[X,[X,Y]]]]]]\\
      - 6618*[[X,Y],[Y,[X,[X,[X,[X,Y]]]]]] - 4940*[[X,Y],[Y,[Y,[X,[X,[X,Y]]]]]]\\
      - 226*[[X,Y],[Y,[Y,[Y,[X,[X,Y]]]]]] + 264*[[X,Y],[Y,[Y,[Y,[Y,[X,Y]]]]]]\\
      - 5712*[[X,Y],[[X,Y],[X,[X,[X,Y]]]]] - 6480*[[X,Y],[[X,Y],[Y,[X,[X,Y]]]]]\\
      - 1188*[[X,Y],[[X,Y],[Y,[Y,[X,Y]]]]] - 2160*[[X,[X,Y]],[X,[X,[X,[X,Y]]]]]\\
      - 3852*[[X,[X,Y]],[Y,[X,[X,[X,Y]]]]] - 1970*[[X,[X,Y]],[Y,[Y,[X,[X,Y]]]]]\\
      - 1166*[[X,[X,Y]],[Y,[Y,[Y,[X,Y]]]]] - 1464*[[X,[X,Y]],[[X,Y],[X,[X,Y]]]]\\
      + 3280*[[X,[X,Y]],[[X,Y],[Y,[X,Y]]]] - 3738*[[Y,[X,Y]],[X,[X,[X,[X,Y]]]]]\\
      - 3360*[[Y,[X,Y]],[Y,[X,[X,[X,Y]]]]] + 2356*[[Y,[X,Y]],[Y,[Y,[X,[X,Y]]]]]\\
      + 750*[[Y,[X,Y]],[Y,[Y,[Y,[X,Y]]]]] - 5520*[[Y,[X,Y]],[[X,Y],[X,[X,Y]]]]\\
      - 356*[[Y,[X,Y]],[[X,Y],[Y,[X,Y]]]] + 2412*[[X,[X,[X,Y]]],[Y,[X,[X,Y]]]]\\
      + 580*[[X,[X,[X,Y]]],[Y,[Y,[X,Y]]]] - 2884*[[Y,[X,[X,Y]]],[Y,[Y,[X,Y]]]]
\end{multline}
\normalsize

Considérons la base des mots cycliques de longueur 8 suivante:\\

\verb"C8[1] := [1, 1, 1, 1, 1, 1, 1, 1]; C8[2] := [1, 1, 1, 1, 1, 1, 1, 2];"

 \verb"C8[3] := [1, 1, 1, 2, 1, 1, 1, 2]; C8[4] := [1, 1, 1, 1, 2, 1, 1, 2];"

  \verb"C8[5] := [1, 1, 1, 1, 1, 2, 1, 2]; C8[6] := [1, 1, 2, 1, 1, 2, 1, 2];"

   \verb"C8[7] := [1, 1, 1, 2, 1, 2, 1, 2]; C8[8] := [1, 2, 1, 2, 1, 2, 1, 2];"

    \verb"C8[9] := [1, 1, 1, 1, 1, 1, 2, 2]; C8[10] := [1, 1, 1, 1, 2, 2, 1, 2];"

     \verb"C8[11] := [1, 1, 1, 2, 2, 1, 1, 2]; C8[12] := [1, 1, 1, 2, 1, 1, 2, 2];"

      \verb"C8[13] := [1, 1, 2, 2, 1, 2, 1, 2]; C8[14] := [1, 1, 2, 2, 1, 1, 2, 2];"

       \verb"C8[15] := [1, 1, 1, 1, 2, 1, 2, 2]; C8[16] := [1, 1, 2, 1, 2, 2, 1, 2];"

        \verb"C8[17] := [1, 1, 2, 1, 2, 1, 2, 2]; C8[18] := [1, 1, 1, 2, 2, 1, 2, 2];"

         \verb"C8[19] := [1, 2, 1, 2, 2, 1, 2, 2]; C8[20] := [1, 1, 1, 1, 1, 2, 2, 2];"

          \verb"C8[21] := [1, 1, 1, 2, 2, 2, 1, 2]; C8[22] := [1, 1, 2, 1, 1, 2, 2, 2];"

           \verb"C8[23] := [1, 1, 2, 2, 2, 1, 2, 2]; C8[24] := [1, 1, 1, 2, 1, 2, 2, 2];"

            \verb"C8[25] := [1, 2, 1, 2, 1, 2, 2, 2]; C8[26] := [1, 1, 2, 2, 1, 2, 2, 2];"

             \verb"C8[27] := [1, 2, 2, 2, 1, 2, 2, 2]; C8[28] := [1, 1, 1, 1, 2, 2, 2, 2];"

              \verb"C8[29] := [1, 1, 2, 2, 2, 2, 1, 2]; C8[30] := [1, 1, 2, 1, 2, 2, 2, 2];"

               \verb"C8[31] := [1, 2, 2, 1, 2, 2, 2, 2]; C8[32] := [1, 1, 1, 2, 2, 2, 2, 2];"

                \verb"C8[33] := [1, 2, 1, 2, 2, 2, 2, 2]; C8[34] := [1, 1, 2, 2, 2, 2, 2, 2];"

                 \verb"C8[35] := [1, 2, 2, 2, 2, 2, 2, 2]; C8[36] := [2, 2, 2, 2, 2, 2, 2, 2]"

\vspace{0,3cm}

Le terme d'ordre 8 du membre de droite de~(\ref{KV2}) dans la base cyclique ci-dessus vaut exactement d'après le programme Maple (au facteur $\frac1{21772800}$ près)

\begin{multline}
[0,72,720,-1080,864,-2160,4320,4320,-252,\textbf{-1080,0,0},\\-2160,540,\textbf{-1080},-2160,-2160,
0,-2160,504,1440,0,\textbf{0},1440,4320,\\\textbf{0},720,-630,\textbf{-1080,-1080},-1080,504,864,-252,72,0]
\end{multline}

Le programme Caml a calculé, pour la composante de degré $8$ du  couple $(F^0, G^0)$, les coordonnées dans la base cyclique ci-dessus, de $$\ad X\circ \partial_X F+\ad Y\circ \partial_Y G$$et trouve (au facteur $ \frac1{21772800}$ près)
\begin{multline}
[0, 72, 720, -1080, 864, -2160, 4320, 4320, -252, \textbf{-1088, 40, -40},\\ -2160,
    540, \textbf{-1072}, -2160, -2160, 0, -2160, 504, 1440, 0, \textbf{-40}, 1440, 4320,\\ \textbf{40},
    720, -630, \textbf{-1072, -1088}, -1080, 504, 864, -252, 72, 0]
\end{multline}

Comme on le constate, la différence est non nulle et vaut (au facteur $\frac1{21772800}$ près)

\begin{eqnarray} \label{diff}
8*[1,1,1,1,2,2,1,2] -8*[1,1,1,1,2,1,2,2] -8*[1, 1, 2, 2, 2, 2, 1, 2]\\\nonumber + 8*[1,1,2,1,2,2,2,2]
+ 40*[1,1,1,2,1,1,2,2]-40*[1,1,1,2,2,1,1,2]\\\nonumber +40*[1,1,2,2,2,1,2,2]-40*[1,1,2,2,1,2,2,2],
\end{eqnarray}

\noindent c'est à dire en notant $x=\ad X, y=\ad Y$,
\begin{equation}\nonumber 8*\tr_\g\left(x^4y^2xy-x^4yxy^2 -x^2y^4xy +x^2yxy^4 \right) + 40* \tr_\g\left(x^3yx^2y^2-x^3y^2x^2y+x^2y^3xy^2-x^2y^2xy^3 \right)
\end{equation}

On en déduit le résultat annoncé:

\begin{theo} La solution proposée par Kashiwara et Vergne \cite{KV} dans le cas résoluble à la conjecture KV-78 ne convient pas en degré $8$ pour les algèbres de Lie générales.
\end{theo}

\noindent \textit{Preuve :} Notons $x=\ad X, y=\ad Y$. Considérons les termes de degré $3$ en $y$. Il suffit de trouver une algèbre de Lie telle que $\tr_\g( x^4y^2xy - x^4yxy^2 ) + 5* \tr_\g\left(x^3yx^2y^2-x^3y^2x^2y\right) \neq 0
$. Considérons $\g=\mathrm{gl}(n, \mathbf{R})\times V$ avec $V$ une représentation telle que $V^*$ ne soit pas isomorphe à~$V$. On prendra pour $V$ la représentation standard.

Pour $X,Y$ deux matrices, l'action adjointe dans $\g$ laisse stable $\mathrm{gl}(n,\mathbf{R})$ et $V$. La composante sur $\mathrm{gl}(n,\mathbf{R})$ de la trace sera nulle car $\mathrm{gl}(n,\mathbf{R})$ est une algèbre quadratique, et donc nos expressions valent les traces matricielles standards. On peut prendre par exemple

$$X:=\left(
  \begin{array}{ccc}
    1 & 1 & 0  \\
    0 & 0 & 1  \\
    0 & 0 & 0  \\
  \end{array}
\right)
\quad
Y:=  \left(
  \begin{array}{ccccc}
    0 & 0 & 0  \\
    1 & 0 & 0  \\
    0 & -1 & 0  \\
  \end{array}
\right).$$

On trouve $\tr(X^4Y^2XY)=-1$,  $\tr(X^4YXY^2)=1$, $\tr(X^3YX^2Y^2) = \tr(X^3Y^2X^2Y)=-1$ ce qui permet de conclure.

\hfill $\blacksquare$

\paragraph{Remarque dans le cas quadratique :} La différence~(\ref{diff}) est nulle  pour les algèbres quadratiques. On a vérifié  que les termes d'ordre 10 du membre de gauche et du membre de droite  (\ref{KV2}) sont distincts, mais que leur différence s'annule dans le cas quadratique. On peut donc penser que le couple \textit{résoluble} de Kashiwara-Vergne, résout aussi les équations  dans le cas des algèbres quadratiques. Ceci a été vérifié pour  $\mathrm{sl}(2,\mathbf{R})$ dans \cite{rou81}. Dans le cas quadratique \cite{Ve} et \cite{AM2} proposent des solutions  distinctes de la solution résoluble. Ces deux solutions du cas quadratique ne sont pas universelles d'après \cite{AP}.

\end{document}